\documentclass[10pt, a4paper]{amsart}

\usepackage{amsmath}
\usepackage{wasysym}
\usepackage{amsthm}
\usepackage{latexsym}
\usepackage{amsfonts}
\usepackage{mathrsfs}
\usepackage{amssymb}
\usepackage{dsfont}
\usepackage[all]{xy}
\usepackage{color}
\usepackage{cleveref}

\newtheorem{Theorems1}{Theorem}[section]

\newtheorem{Coroll1}[Theorems1]{Corollary}
\newtheorem{Lemma1}[Theorems1]{Lemma}
\newtheorem{Proposition1}[Theorems1]{Proposition}

\theoremstyle{definition}
\newtheorem{Definitions1}[Theorems1]{Definition}

\theoremstyle{plain}

\theoremstyle{remark}

\theoremstyle{plain}

\newtheorem*{Ques}{Question}

\newtheorem*{Theo}{Theorem}

\newtheorem*{Cor}{Corollary}

\theoremstyle{remark}

\theoremstyle{plain}

\theoremstyle{remark}
\newtheorem*{Rem}{Remark}

\numberwithin{equation}{Theorems1}

\numberwithin{equation}{Theorems1}
\def\multiset#1#2{\ensuremath{\left(\kern-.3em\left(\genfrac{}{}{0pt}{}{#1}{#2}\right)\kern-.3em\right)}}

\def\Multiset#1#2{\left(\!\!\left(\!\!
    \begin{array}{c}
      n \\
      r
    \end{array}
  \!\!\right)\!\!\right)}

\newcommand{\spec}[1]{\operatorname{Spec}(#1)}
\newcommand{\pt}[1]{\operatorname{pt}(#1)}


\begin{document}

\title{A homological reformulation of the link condition}
\author{Rishi Vyas}
\address{Wolfson College, Cambridge, CB3 9BB, United Kingdom}
\date{}

\begin{abstract}
We prove an equivalent condition for the existence of a link between prime ideals in terms of the structure of a certain cohomology module. We use this formulation to answer an open question regarding the nature of module extensions over one sided noetherian rings. We apply the techniques developed in this paper to the local link structure of prime ideals of small homological height and examine when certain noetherian rings satisfy the density condition.
\end{abstract}


\maketitle

All rings in this paper are associative with a unit.

When can one localize at a prime ideal in a noncommutative noetherian ring? It is well known that this is not always possible. The key insight of M\"{u}ller, Jategaonkar and others was to realize that, in some situations, the obstructions to localization are other prime ideals and other prime ideals alone. For an overview of this theory, we recommend \cite{GW} or \cite{MR}; another good reference is Jategaonkar's monograph, \cite{JAT}. There exists a relation on $\spec{R}$, called a \textit{link} (denoted by $\rightsquigarrow$), which measures this obstruction; if $X$ is a right Ore set contained in $C_{R}(Q)$, the set of elements of $R$ regular modulo $Q$, and $P$ is a prime ideal such that $P\rightsquigarrow Q$, then $X$ is contained in $C_{R}(P)$ \cite[Lemma 14.17]{GW}. Determining this relation, however, can be a reasonably involved task. The existence of these links can be deduced from certain module extensions, and in order to ensure that suitable module extensions allow us to detect links, we must ensure that our rings satisfy the second layer condition. Given a module over a ring which satisfies the second layer condition, it is possible to break $M$ up into layers, each of which is a fully faithful module over some prime factor of $R$, such that the associated primes are linked to each other. The reader can find more details about the second layer condition in \cite{GW} or \cite{JAT}.

In this paper, we attempt to give a homological criterion for this relation.  Much of what is done below is just a translation of what is known into a different language, but in a setting where the tools of homological algebra are then available to work with. This approach leads to certain applications. It is well known that uniform modules over right noetherian rings come in two distinct flavours: tame and wild. Jategaonkar's main lemma, \cite[Theorem 12.1]{GW}, can be interpreted as saying that if we take a tame uniform module $U$ over a noetherian ring and quotient by the affiliated submodule corresponding to its unique associated prime, in many cases the resulting quotient module (the `second layer' of $U$) contains only tame submodules; we prove the corresponding statement for one sided noetherian rings in \cref{mainlemmafull}. One can ask the corresponding question for wild modules: is the second layer of a wild uniform module over a noetherian ring wild? We show that this is indeed the case for many families of noetherian rings. This gives information about the structure of indecomposable injective modules over noetherian rings. Another application is in \textsection 5, where we describe the link structure of `small' prime ideals in certain rings. 

In \textsection 1, we define the truncated local cohomology of a module at an ideal, and show how the truncated local cohomology at the `point' of an injective module contains all possible second layers. This `global' approach is then used to give a new proof of the fact that second layers of torsionfree modules are torsionfree. We prove a version of Jategaonker's main lemma, \cite[Theorem 12.1]{GW}, as  \cref{mainlemmafull}:

\begin{Theo}
Let $R$ be a right noetherian ring, and let $M$ be an $R$-module with an affiliated series $0\lneq U\lneq_{e} M$, with corresponding affiliated primes $Q$ and $P$. Then, one of the following two situations occur, and they are mutually exclusive:
\begin{itemize}
\item[i.)] $P\lneq Q$, and $U\lneq \mathrm{ann}_{M}(P)$.
\item[ii.)] $M/U$ embeds into $\mathrm{Ext}^{1}_{R}(R/P,\mathrm{ann}_{E}(Q))$, where $E$ is the injective hull of $M$. There is an indecomposable injective module $E'$ with assassinator $Q$ such that $\mathrm{Ext}^{1}_{R}(R/P,\mathrm{ann}_{E'}(Q))$ contains a fully faithful $R/P$-module. In this case, if $U$ is a torsionfree $R/Q$-module, then $M/U$ is a torsionfree $R/P$-module.
\end{itemize}
\end{Theo}

Unlike the classical arguments, our arguments work in the context of one sided noetherian rings. This allows us to answer the open question posed at the end of \cite[Exercise 12P, pg. 219]{GW} in the affirmative. 

In \textsection 2 we reprove some known results from a slightly different perspective. While the methods in this section were arrived at independently by the author, some of the ideas used are hinted at by Brown and Warfield in \cite{KW} and by Brown in \cite{KB2}. The point of this section is to show that it is possible to deduce the existence of a link by directly considering the module theoretic structure of our truncated local cohomology modules. 

In \textsection 3, we study how the structure of our local cohomology modules impacts the link structure of our clique. It is a major and long-standing open question in the theory of noetherian rings as to whether links can form between comparable prime ideals. We establish some conditions for incomparability in terms of the structure of our local cohomology modules. The following appears as \cref{INC}:

\begin{Theo}
Let $R$ be a right noetherian ring, and let $Q$ be a prime ideal. Then, the following are equivalent:
\begin{itemize}
\item[i.)] There does not exist a prime ideal $P$ such that $P\gneq Q$ and $P\rightsquigarrow Q$.
\item[ii.)]  The set of left torsion elements of $Q/Q^{2}$ is contained in the right torsion subbimodule of $Q/Q^{2}$.
\item[iii.)] $\mathrm{Ext}^{1}_{R}(R/Q,\mathcal{Q}(R/Q))$ is injective as a right $R/Q$-module.
\item[iv.)] $\mathrm{Tor}_{1}^{R/Q}(R/I,Q/Q^{2})$ is torsion as a right $R/Q$-module for any left ideal $I\leq R/Q$.
\item[v.)] $\mathrm{Ext}^{1}_{R}(R/P,\mathcal{Q}(R/Q))=0$ for all $P\gneq Q$.
\item[vi.)] $\mathrm{Ext}^{1}_{R}(\text{---},\mathcal{Q}(R/Q))$ is an exact functor on the category of $R/Q$-modules.
\end{itemize}
\end{Theo}

Some of the arguments that we use to establish the above result can also be recycled in our attempt to study how the structure of the tame indecomposable injective module associated to a prime ideal $P$ impacts the structure of the various wild indecomposable injectives associated to $P$. Work has been done in the past on converses to Jategaonkar's main lemma for wild modules (see \cite{KW}). We prove some results of a similar flavour, for example:

\begin{Cor}[\cref{impact2}]
Let $R$ be a right noetherian ring. If $\mathrm{Ext}^{1}_{R}(R/P,\mathcal{Q}(R/Q))\neq 0$, then $\mathrm{Ext}^{1}_{R}(R/P,\mathrm{ann}_{E}(Q))\neq 0$ for all indecomposable injective modules $E$ with assassinator $Q$.
\end{Cor}

In \textsection 4, we study the links that arise between small prime ideals for suitably nice rings. We show that for a reasonably large class of rings, the structure of a clique of prime ideals of homological height $1$ is very simple. We also prove a result regarding the structure of a clique of prime ideals of homological height $2$. In particular, we show that there is a map on the spectrum of our ring $R$, often induced by a distinguished automorphism of $R$, which plays a role in determining the localizability of small prime ideals.

 Jategaonkar showed, in \cite[Theorem 3.5]{JAT2}, that the second layers of wild modules are wild if and only if $R$ happened to satisfy the \textit{right density condition}. In \textsection 5, we see how homological restrictions on our ring allow us to study questions regarding the density of certain modules. We show that a large number of rings satisfy the right density condition. Using \cref{listofloc}, \cref{rdc3} in particular proves the following:

\begin{Theo}
Let $k$ be a field. Let $A$ be a $k$-algebra such that either:
\begin{itemize}
\item [a.)] $A$ has a noetherian connected filtration such that $gr(A)$ is commutative of finite global dimension.
\item [b.)] $A$ is a complete, Auslander-Regular, local algebra with maximal ideal $m$ such that $A/m$ is finite dimensional over $k$.
\end{itemize}
Suppose every clique in $A$ is localizable. Then $A$ satisfies the right density condition. In particular, the second layer of a wild indecomposable injective module over $A$ is wild.
\end{Theo}

Unless we explicitly mention otherwise, modules over $R$ will be right modules. We will always assume that $R$ is right noetherian, though we continue to mention it in the hypotheses of our results for the sake of clarity, and we will explicitly mention the fact that $R$ is noetherian (by which we mean right and left noetherian) whenever we require it.

\section{Truncated Local Cohomology} \label{TLCCHAP}

Any uniform module $U$ over a right noetherian ring $R$ has a unique associated prime, which is then called the assassinator of $U$. Note that an injective module is uniform if and only if it is indecomposable. We remind the reader about some notation: if $M$ is an $R$-module and $I$ an ideal of $R$, $\mathrm{ann}_{M}(I):=\{x\in M| xI=0\}$.

\begin{Definitions1} \label{point}
Let $R$ be a ring. Let $U$ be a uniform module over $R$, with assassinator $Q$. Define $\operatorname{pt}(U):= \mathrm{ann}_{U}(Q)$. We will call $\operatorname{pt}(U)$ the `point' of $U$.
\end{Definitions1}

\begin{Definitions1} \label{secondlayerdef}
Let $R$ be a ring. Let $U$ be a uniform module over $R$, with assassinator $Q$. The `second layer' of $U$ is the module $U/\operatorname{pt}(U)$.
\end{Definitions1}

The reader should be aware that there are other definitions of the `second layer' of a module in the literature which do not coincide with ours. Indeed, our definition is not a very good one --- a better definition of the second layer of a uniform module $U$ would involve the affiliated primes of $U/\operatorname{pt}(U)$. Our definition is purely for bookkeeping purposes; we will not need anything more complicated in this paper.

It is immediate that for any $R$-module $M$, there is a natural isomorphism of $R/I$-modules $\mathrm{Hom}_{R}(R/I,M)\cong \mathrm{ann}_{M}(I)$ for any ideal $I$. Thus, if $E$ is an indecomposable injective module over $R$ with assassinator $P$, $\pt{E}$ is an indecomposable injective over $R/P$. 

Given a prime ideal $P$ in $R$, let $E_{P}$ be the injective hull of any uniform right ideal of $R/P$ over $R$. It can be shown that this construction, up to isomorphism, is independent of the choice of right ideal, and $E_{P}$ can be characterized as the unique indecomposable injective with assassinator $P$ such that $\pt{E_{P}}$ is torsionfree as an $R/P$-module. We can explicitly describe $\pt{E_{P}}$: it is the unique simple module of the Goldie quotient ring of $R/P$. 

In general, an indecomposable injective module is called \textit{tame} if it is $E_{Q}$ for some prime $Q$, and \textit{wild} otherwise. We extend this definition by saying that a uniform module is tame (resp. wild) if its injective hull is tame (resp. wild). 

\begin{Definitions1} \label{localcoho}
Let $R$ be a ring. Let $I$ be an ideal of $R$. The ${i^{th}\ truncated\  local\  cohomology}$ of $M$ at $I$, $\mathrm{\mathrm{T}\Gamma}_{I}^{i}(M)$, is $\mathrm{Ext}^{i}_{R}(R/I,M)$. 
\end{Definitions1}

By abuse of language, we will usually refer to $\mathrm{\mathrm{T}\Gamma}_{I}^{1}(M)$ as the local cohomology of $M$ at $I$ when there is no chance of confusion. When $P$ is a prime ideal,  $\mathrm{\mathrm{T}\Gamma}_{P}^{i}(M)$ is a useful invariant of $M$. Its structure as an $R/P$-module can be used to determine the multiplicity of $E_{P}$ in the $i^{th}$ term of a minimal injective resolution of $M$ --- to be precise, we have the following result:

\begin{Lemma1}\label{injnum}
 Let $R$ be a right noetherian ring. Let $M$ be a right $R$-module with minimal injective resolution $I$. Then, the multiplicity of $E_{P}$ as a summand of $I^{i}$ is the reduced rank of $\mathrm{Ext}^{i}_{R}(R/P,M)$ as a right $R/P$-module. In particular, $E_{P}$ appears as a summand in $I^{i}$ if and only if $\mathrm{Ext}^{i}_{R}(R/P,M)$ is not torsion as an $R/P$-module.
\begin{proof}
\cite[Lemma 2.5]{JY}.
\end{proof}
\end{Lemma1}

As $\mathrm{\mathrm{T}\Gamma}_{I}$ is a $\delta$-functor, we have, for any exact sequence of modules $0\rightarrow U\rightarrow V\rightarrow V/U\rightarrow 0$ a corresponding long exact sequence of $R/I$-modules

\begin{equation}\label{exact}
0\rightarrow \mathrm{\mathrm{T}\Gamma}_{I}^{0}(U)\rightarrow \mathrm{\mathrm{T}\Gamma}_{I}^{0}(V)\rightarrow \mathrm{\mathrm{T}\Gamma}_{I}^{0}(V/U)\rightarrow \mathrm{\mathrm{T}\Gamma}_{I}^{1}(U)\rightarrow \mathrm{\mathrm{T}\Gamma}_{I}^{1}(V)\rightarrow ...
\end{equation}

We will use the notation $\partial$ to denote the boundary map in the above long exact sequence. 

Let us describe the behaviour of this exact sequence in a special situation. We will begin by giving two versions of a general lemma showing the annihilators of the $n^{th}$ `layer' of a module at a prime ideal $P$ sit inside suitable  $\mathrm{Ext}$ modules. The first variant is essentially a computation of local cohomology at a power of a prime ideal $P$. The second is just a restatement of the first version, but in a form where the curious dichotomy of the situation becomes apparent; this is the form that we will use for our homological formulation of the link condition. 

\begin{Lemma1} \label{primmainlemmam}
 Let $R$ be a right noetherian ring, and suppose $V$ is an $R$-module with a series $0\lneq U\lneq_{e} V$, with $U$ the affiliated submodule corresponding to the affiliated prime $Q$. Let $P$ be a prime ideal in $R$. Then, one of the following two situations occur, and they are mutually exclusive:
\begin{itemize}
\item [i.)] $P\nleq Q$, and the map induced by $\partial$ embeds $\mathrm{\mathrm{T}\Gamma}_{P^{n}Q}^{0}(V)/\mathrm{\mathrm{T}\Gamma}_{Q}^{0}(V)$ into $\mathrm{\mathrm{T}\Gamma}_{P^{n}}^{1}(U)$ for all $n$.
\item [ii.)] $P\leq Q$, and the map induced by $\partial$ embeds $\mathrm{\mathrm{T}\Gamma}_{P^{n}Q}^{0}(V)/\mathrm{\mathrm{T}\Gamma}_{P^{n}}^{0}(V)$ into $\mathrm{\mathrm{T}\Gamma}_{P^{n}}^{1}(U)$ for all $n$.
\end{itemize}
\begin{proof}
This follows from examining the long exact sequence, as in \cref{exact}, arising from the short exact sequence $0\rightarrow U\rightarrow V\rightarrow V/U\rightarrow 0$ with respect to the ideal $P^{n}$. Recall that the map $\partial$ in the statement of the lemma is the boundary map in the above long exact sequence.

First, note that $\mathrm{\mathrm{T}\Gamma}_{P^{n}}^{0}(V/U) \cong \mathrm{\mathrm{T}\Gamma}_{P^{n}Q}^{0}(V)/U$. If $P\nleq Q$, then $\mathrm{T\Gamma}_{P^{n}}^{0}(U)$ and thus $\mathrm{T\Gamma}_{P^{n}}^{0}(V)$ are both $0$. Note that $U\cong \mathrm{T\Gamma}_{Q}^{0}(V)$, as it is the affiliated submodule corresponding to the affiliated prime $Q$. Thus, the first situation holds.

However, if $P\leq Q$, then by examining the long exact sequence and arguing in a similar fashion, we see that the second situation holds.
\end{proof}
\end{Lemma1}

This lemma allows us to see the analogy with Matlis' theory of indecomposable injectives over a commutative noetherian ring (see \cite[pg. 99]{LAM}). We indulge ourselves by going through the details.  

\begin{Coroll1} \label{matlisfornoncomm}
Let $R$ be a right noetherian ring, and let $P$ be a prime ideal. Let $E$ be an indecomposable injective module with associated prime $Q$. Then, there are the following isomorphisms:
\begin{itemize}
\item[i.)] If $P\nleq Q$, then $\mathrm{\mathrm{T}\Gamma}_{P^{n}Q}^{0}(E)/\pt{E}\cong \mathrm{\mathrm{T}\Gamma}_{P^{n}}^{1}(\pt{E})$ for all $n$.
\item[ii.)] If $P\leq Q$,  then $\mathrm{\mathrm{T}\Gamma}_{P^{n}Q}^{0}(E)/\mathrm{\mathrm{T}\Gamma}_{P^{n}}^{0}(E)\cong \mathrm{\mathrm{T}\Gamma}_{P^{n}}^{1}(\pt{E})$ for all $n$.
\end{itemize}
\begin{proof}
The fact that $E$ is injective tells us that the maps in \cref{primmainlemmam} are isomorphisms. The result follows.
\end{proof}
\end{Coroll1}

For the rest of this paper, we will only be concerned with local cohomology at a prime ideal $P$. 

\begin{Lemma1} \label{primmainlemma}
 Let $R$ be a right noetherian ring, and suppose $V$ is an $R$-module with a series $0\lneq U\lneq_{e} V$, with $U$ the affiliated submodule corresponding to the affiliated prime $Q$. Let $P$ be a prime ideal in $R$. Then, one of the following two situations occur, and they are mutually exclusive:
\begin{itemize}
\item [i.)] $P\lneq Q$, and $\mathrm{\mathrm{T}\Gamma}_{P}^{0}(U)\lneq \mathrm{\mathrm{T}\Gamma}_{P}^{0}(V)$.
\item [ii.)]The boundary map $\partial:\mathrm{\mathrm{T}\Gamma}_{P}^{0}(V/U)\rightarrow \mathrm{\mathrm{T}\Gamma}_{P}^{1}(U)$ is injective.
\end{itemize}
\begin{proof}
The boundary map is injective if and only if $\mathrm{\mathrm{T}\Gamma}_{P}^{0}(U)= \mathrm{\mathrm{T}\Gamma}_{P}^{0}(V)$. If $P$ is not contained in $Q$, $\mathrm{\mathrm{T}\Gamma}_{P}^{0}(U)$ and thus $\mathrm{\mathrm{T}\Gamma}_{P}^{0}(V)$ are $0$ as $U$ is a fully faithful $R/Q$-module and $U\lneq_{e}V$. If $P=Q$, then $\mathrm{\mathrm{T}\Gamma}_{P}^{0}(U) = \mathrm{\mathrm{T}\Gamma}_{P}^{0}(V)$. In both cases, the long exact sequence informs us that $\partial:\mathrm{\mathrm{T}\Gamma}_{P}^{0}(V/U)\rightarrow \mathrm{\mathrm{T}\Gamma}_{P}^{1}(U)$ is injective. It follows that if the second situation does not  occur, then $P\lneq Q$ and $\mathrm{\mathrm{T}\Gamma}_{P}^{0}(U)\lneq \mathrm{\mathrm{T}\Gamma}_{P}^{0}(V)$.
\end{proof}
\end{Lemma1}

Matters become even sharper when we restrict ourselves to injective modules.

\begin{Proposition1} \label{mainlemma}
 Let $R$ be a right noetherian ring, and consider an indecomposable injective $E$ with associated prime $Q$. Let $P$ be a prime ideal in $R$. Then, one of the following two situations occur, and they are mutually exclusive:
\begin{itemize}
 \item[i.)] $P\lneq Q$, and there is an exact sequence $$0\rightarrow \mathrm{\mathrm{T}\Gamma}_{P}^{0}(E)/\pt{E}\rightarrow \mathrm{\mathrm{T}\Gamma}_{P}^{0}(E/\pt{E})\rightarrow \mathrm{\mathrm{T}\Gamma}_{P}^{1}(\pt{E})\rightarrow 0.$$
 \item[ii.)] The boundary map $\partial: \mathrm{\mathrm{T}\Gamma}_{P}^{0}(E/\pt{E})\rightarrow 
\mathrm{\mathrm{T}\Gamma}_{P}^{1}(\pt{E})$ is an isomorphism.
\end{itemize}
 If $E=E_{Q}$, then the torsion submodule of $\mathrm{\mathrm{T}\Gamma}_{P}^{1}(\pt{E_{Q}})$ is  unfaithful as an $R/P$-module.

\begin{proof}
 Note that $\mathrm{\mathrm{T}\Gamma}_{P}^{1}(E)=0$ as $E$ is an injective $R$-module. Thus the boundary map is surjective. If $P$ is not contained in $Q$ or is equal to $Q$, the argument follows as in \cref{primmainlemma} and we see that the boundary map is an isomorphism. If $P\lneq Q$, suppose that $\mathrm{\mathrm{T}\Gamma}_{P}^{0}(\pt{E}) = \mathrm{\mathrm{T}\Gamma}_{P}^{0}(E)$. As $P\lneq Q$, $\mathrm{\mathrm{T}\Gamma}_{P}^{0}(\pt{E}) = \pt{E}$. Now, $\mathrm{\mathrm{T}\Gamma}_{P}^{0}(E)$ is an injective $R/P$-module, and thus $\pt{E}$ is an injective, unfaithful $R/P$-module. However, over a prime right noetherian ring, injective modules are faithful \cite[Lemma 14.9]{GW}. Thus $\mathrm{\mathrm{T}\Gamma}_{P}^{0}(\pt{E}) \lneq \mathrm{\mathrm{T}\Gamma}_{P}^{0}(E)$, and the boundary map is not injective. This completes the argument for the first part of the proposition.

Now, let $E=E_{Q}$. Since $\pt{E_{Q}}$ is the unique simple module for $\mathcal{Q}(R/Q)$, the Goldie quotient ring of $R/Q$, $\mathrm{\mathrm{T}\Gamma}_{P}^{1}(\pt{E_{Q}})$ is a summand of $\mathrm{\mathrm{T}\Gamma}_{P}^{1}(\mathcal{Q}(R/Q))$. This is a $(\mathcal{Q}(R/Q),R/P)$-bimodule, and is finitely generated on the left. By \cite[Lemma 8.3]{GW}, the torsion submodule of $\mathrm{\mathrm{T}\Gamma}_{P}^{1}(\mathcal{Q}(R/Q))$ is unfaithful, and the same therefore holds for $\mathrm{\mathrm{T}\Gamma}_{P}^{1}(\pt{E_{Q}})$. This proves the second assertion. 
\end{proof}
\end{Proposition1}

The following corollaries are an immediate consequence of the above proposition. 

\begin{Coroll1} \label{ff1}
Let $R$ be a right noetherian ring, and suppose that $V$ is a uniform module with an affiliated series $0\lneq U\lneq_{e} V$, with corresponding affiliated primes $Q$ and $P$. If $\mathrm{ann}_{U}(P)= \mathrm{ann}_{V}(P)$, then $V/U$ embeds into $\mathrm{\mathrm{T}\Gamma}_{P}^{1}(\pt{E})$, where $E$ is the injective hull of $V$. In particular, $\mathrm{\mathrm{T}\Gamma}_{P}^{1}(\pt{E})$ contains a fully faithful $R/P$-submodule.
\begin{proof}
We know that $V$ embeds into some indecomposable injective module $E$ with assassinator $Q$, and thus that $V/U$ embeds into $E/\pt{E}$. If $P$ is not properly contained in $Q$, then the second situation in \cref{mainlemma} occurs, and $V/U$ embeds in $\mathrm{\mathrm{T}\Gamma}_{P}^{1}(\pt{E})$. Thus $\mathrm{\mathrm{T}\Gamma}_{P}^{1}(\pt{E})$ contains a fully faithful submodule. If, on the other hand, $P\lneq Q$, then we are in the first situation in \cref{mainlemma}. There is an exact sequence $$0\rightarrow \mathrm{ann}_{E}(P)/\pt{E}\rightarrow \mathrm{\mathrm{T}\Gamma}_{P}^{0}(E/\pt{E})\rightarrow \mathrm{\mathrm{T}\Gamma}_{P}^{1}(\pt{E})\rightarrow 0.$$ Our conditions on $V$ imply that $V\cap \mathrm{ann}_{E}(P) = U$. From the above exact sequence, we again see that $V/U$ embeds in $\mathrm{\mathrm{T}\Gamma}_{P}^{1}(\pt{E})$.
\end{proof}
\end{Coroll1} 

\begin{Coroll1} \label{tf}
Let $R$ be a right noetherian ring, and suppose that $V$ is a uniform module with an affiliated series $0\lneq U\lneq_{e} V$, with corresponding affiliated primes $Q$ and $P$. Suppose $U$ is a torsionfree $R/Q$-module. If $\mathrm{ann}_{U}(P)= \mathrm{ann}_{V}(P)$, then $V/U$ is a torsionfree $R/P$-module. 
\begin{proof}
As in \cref{ff1}, $V/U$ embeds into $\mathrm{\mathrm{T}\Gamma}_{P}^{1}(\pt{E_{Q}})$. By \cref{mainlemma}, the torsion submodule of $\mathrm{\mathrm{T}\Gamma}_{P}^{1}(\pt{E_{Q}})$ is unfaithful. As $V/U$ is fully faithful as an $R/P$-module, it follows that $V/U$ is torsionfree.
\end{proof}
\end{Coroll1}

We note that the above two corollaries hold without the assumption that $V$ is uniform. See \cref{mainlemmafull}.

There is a converse to \cref{tf}. Compare the following with \cite[Theorem 12.2]{GW}.

\begin{Proposition1} \label{conversetf}
Let $R$ be a right noetherian ring, and let $P$ and $Q$ be prime ideals. If $\mathrm{\mathrm{T}\Gamma}_{P}^{1}(\pt{E_{Q}})$ has a torsionfree submodule as an $R/P$-module, then there exists a finitely generated uniform right module $V$ with an affiliated series $0\lneq U\lneq V$ with corresponding affiliated primes $Q$ and $P$ such that $U$ is a finitely generated torsionfree uniform $R/Q$-module and $V/U$ is isomorphic to a uniform right ideal of $R/P$. If $R/Q$ is a left Goldie ring, then $U$ is also isomorphic to a uniform right ideal of $R/Q$.
\begin{proof}
 There are two cases to deal with. Let us first suppose that $P$ is not contained properly in $Q$. Then \cref{mainlemma} tells us that $\partial: \mathrm{\mathrm{T}\Gamma}_{P}^{0}(E_{Q}/\pt{E_{Q}})\rightarrow 
\mathrm{\mathrm{T}\Gamma}_{P}^{1}(\pt{E_{Q}})$ is an isomorphism. As $\mathrm{\mathrm{T}\Gamma}_{P}^{1}(\pt{E_{Q}})$ is not torsion over $R/P$, there exists a submodule $W\leq \mathrm{\mathrm{T}\Gamma}_{P}^{0}(E_{Q}/\pt{E_{Q}})$ which we can take to be isomorphic to a uniform right ideal of $R/P$. This pulls back under the natural projection map to give us a submodule $V'$ of $E_{Q}$ such that $V'/\pt{E_{Q}}\cong W$. Let $x\in V' - \pt{E_{Q}}$. Then, $V=xR$ has an affiliated series $0\lneq U = V\cap \pt{E_{Q}} \lneq V$ with affiliated primes $Q$ and $P$. 

If, on the other hand, $P \lneq Q$, then there is an exact sequence  
$$0\rightarrow \mathrm{\mathrm{T}\Gamma}_{P}^{0}(E_{Q})/\pt{E_{Q}}\rightarrow \mathrm{\mathrm{T}\Gamma}_{P}^{0}(E_{Q}/\pt{E_{Q}})\rightarrow \mathrm{\mathrm{T}\Gamma}_{P}^{1}(\pt{E_{Q}})\rightarrow 0.$$ As  $\mathrm{\mathrm{T}\Gamma}_{P}^{0}(E_{Q}/\pt{E_{Q}})$ is not torsion, we see that there exists a torsionfree $R/P$-submodule $W\leq \mathrm{\mathrm{T}\Gamma}_{P}^{0}(E_{Q}/\pt{E_{Q}})$ --- this implies that $W\cap  \mathrm{\mathrm{T}\Gamma}_{P}^{0}(E_{Q})/\pt{E_{Q}}= 0$. We can take $W$ to be isomorphic to a uniform right ideal of $R/P$. This pulls back to a submodule $V'$ of $E_{Q}$ such that $V'\cap \mathrm{\mathrm{T}\Gamma}_{P}^{0}(E_{Q}) =\pt{E_{Q}}$. We have that $V'/\pt{E_{Q}}\cong W$. We can now proceed as in the previous case to cut $V'$ down to a finitely generated module.

If $R/Q$ is left Goldie, then any torsionfree uniform finitely generated $R/Q$-module is isomorphic to a right ideal of $R/Q$ (\cite[Corollary 7.20]{GW}).
\end{proof}
\end{Proposition1}

Before continuing, let us take a moment to compare what we have been done above to what was known before.

To give the reader a chance to see the analogy, we state the main classical result as it appears in \cite[Theorem 12.1]{GW}. Before doing so, we need to define the notion of a link between prime ideals (\cite[Chapter 5.3]{JAT}). 

\begin{Definitions1} \label{linknoeth}
Let $P$ and $Q$ be prime ideals in $R$, where $R$ is right noetherian. We say that $P$ is \textit{linked} to $Q$, $P\rightsquigarrow Q$, if there exists an ideal $A$ such that $PQ\leq A \lneq P\cap Q$ with $P\cap Q /A$ fully faithful as a left $R/P$-module and torsionfree as a right $R/Q$-module. In this case, $P\cap Q /A$ is called a \textit{linking bimodule}. There is a corresponding definition for left noetherian rings.
\end{Definitions1}

\cite[Lemma 8.1]{GW} tells us that the above definition reduces to the standard one for prime ideals in a noetherian ring when $R$ is noetherian. 

This relation gives $\mathrm{Spec(R)}$ the structure of a directed graph. The connected component of $P$ is called the \textit{clique} of $P$. The prime ideals that appear in the clique of $P$ are the primary obstruction to localizing at $P$. In general, the existence of non-trivial links from a prime ideal tell us that we are unable to localize at only that prime in any sensible way. In nice situations, though, we are often able to localize at the entire clique and reduce to a more favourable position, where our prime is a (but not necessarily the only!) maximal ideal. A good reference for this material is \cite[Chapters 12 and 14]{GW}, or Jategaonkar's own monograph on localization in noetherian rings, \cite{JAT}.

\begin{Theorems1}[Jategaonkar's main lemma] \label{jml}
Let $R$ be a noetherian ring, and let $M$ be an $R$-module with an affiliated series $0\lneq U\lneq_{e} M$, with corresponding affiliated primes $Q$ and $P$. Suppose that $M'$ is a submodule containing $U$ such that the annihilator of $M'$ is maximal amongst annihilators of submodules properly containing $U$. Let $A=\mathrm{ann}_{R}(M').$  Then, one of the following two situations occur:
\begin{itemize}
\item[i.)] $P\lneq Q$, $M'P=0$, and both $M'$ and $M'/U$ are faithful torsion $R/P$-modules.
\item[ii.)] $P\rightsquigarrow Q$, with $P\cap Q /A$ a linking bimodule. In this case, if $U$ is torsionfree as an $R/Q$-module, then $M'/U$ is torsionfree as an $R/P$-module. 
\end{itemize}
\begin{proof}
\cite[Theorem 12.1]{GW}
\end{proof}
\end{Theorems1}

We view \cref{ff1} as a translation of \cref{jml}. Indeed, if one examines the assumptions on $M'$ in \cref{jml}, then by taking $A= \mathrm{ann}_{R}(M')$ as large as possible, we are assuming that every submodule in $M'$ not contained in $U$ has annihilator $A$. If $P$ is not contained in $Q$, then, as in \cref{mainlemma}, there is no problem. If $P$ is contained in $Q$, then the dichotomy expressed in \cref{jml} is precisely the dichotomy between $\mathrm{ann}_{U}(P)= \mathrm{ann}_{M'}(P)$ and $\mathrm{ann}_{U}(P)\lneq \mathrm{ann}_{M'}(P)$. One should view $M'/U$ appearing in $\mathrm{\mathrm{T}\Gamma}_{P}^{1}(\mathrm{\mathrm{T}\Gamma}_{Q}^{0}(E))$, where $E$ is the injective hull of $M'$ as the favourable eventuality, as we then have that our local cohomology modules have fully faithful submodules. In \textsection 2, we will show that this is equivalent to the existence of a link. We state and prove our variant of \cref{jml} explicitly.

\begin{Theorems1} \label{mainlemmafull}
Let $R$ be a right noetherian ring, and let $M$ be an $R$-module with an affiliated series $0\lneq U\lneq_{e} M$, with corresponding affiliated primes $Q$ and $P$. Then, one of the following two situations occur, and they are mutually exclusive:
\begin{itemize}
\item[i.)] $P\lneq Q$, and $\mathrm{\mathrm{T}\Gamma}_{P}^{0}(U)\lneq \mathrm{\mathrm{T}\Gamma}_{P}^{0}(M)$.
\item[ii.)] $M/U$ embeds into $\mathrm{\mathrm{T}\Gamma}_{P}^{1}(\mathrm{\mathrm{T}\Gamma}_{Q}^{0}(E))$, where $E$ is the injective hull of $M$. There is an indecomposable summand of $E$, $E'$,  with assassinator $Q$ such that $\mathrm{\mathrm{T}\Gamma}_{P}^{1}(\pt{E'})$ contains a fully faithful $R/P$-module. In this case, if $U$ is a torsionfree $R/Q$-module, then $M/U$ is a torsionfree $R/P$-module.
\end{itemize}
\begin{proof}
Most of the work has already been done. The fact that $U\leq_{e} M$ tells us that $Q$ is an affiliated prime of $E$, and that $\mathrm{\mathrm{T}\Gamma}_{Q}^{0}(E)\leq_{e} E$. As $R$ is right noetherian, $E$ is a direct sum of indecomposable injectives, say $\oplus E_{i}$. Each $E_{i}$ has assassinator $Q$, and $\mathrm{\mathrm{T}\Gamma}_{Q}^{0}(E)\cong \oplus \pt{E_{i}}$. 

As in \cref{ff1}, we have that situation ii.) arises precisely when  $M/U$ embeds into $\mathrm{\mathrm{T}\Gamma}_{P}^{1}(\mathrm{\mathrm{T}\Gamma}_{Q}^{0}(E))$. In this case, $\mathrm{\mathrm{T}\Gamma}_{P}^{1}(\mathrm{\mathrm{T}\Gamma}_{Q}^{0}(E))\cong \mathrm{\oplus \mathrm{T}\Gamma}_{P}^{1}(\pt{E_{i}})$ has a fully faithful $R/P$-submodule. We can assume that this submodule is finitely generated and uniform, thereby implying that it sits inside $\mathrm{\mathrm{T}\Gamma}_{P}^{1}(\pt{E_{i}})$ for some $i$. 

If $U$ were torsionfree, then in the above decomposition it follows that $E_{i}\cong E_{Q}$ for all $i$. Thus the torsion submodule of $\mathrm{\mathrm{T}\Gamma}_{P}^{1}(\mathrm{T}\Gamma_{Q}^{0}(E))$ is unfaithful. As $M/U$ is fully faithful, we have that it is torsionfree.  
\end{proof} 
\end{Theorems1}

 \cref{mainlemmafull}, as of now, does not give us a link between $P$ and $Q$.  In \cref{jml}, the proof that $M'/U$ is torsionfree when $U$ is torsionfree  applies only to noetherian rings, as we need to go through the linking bimodule. The argument in \cref{mainlemmafull} recovers the situation, and shows that matters behave similarly for one sided noetherian rings; it should be noted that until now, the question of whether matters behave similarly in the one sided noetherian case was open (\cite[pg. 219]{GW}).

\begin{Rem} \label{remjat}
Jategaonkar, in \cite[pg. 188]{JAT}, offers the following definition of the second layer condition for a prime ideal in a right noetherian ring: a prime ideal $P$ in a right noetherian ring satisfies the right second layer condition if every uniform submodule of the second layer of $E_{P}$ is tame. However, one can also be tempted to define the second layer condition in a manner analogous to how it is defined for noetherian rings, i.e. we can say that a prime ideal $P$ satisfies the right second layer condition if, under the hypothesis of \cref{mainlemmafull}, situation ii.) always occurs when $U$ is torsionfree as an $R/Q$-module.  It follows from \cref{mainlemmafull} that this definition is equivalent to Jategaonkar's. In a similar fashion, we say that a prime ideal $P$ in a right noetherian ring satisfies the strong second layer condition if, under the hypothesis of \cref{mainlemmafull}, situation ii.) always occurs.
\end{Rem}

\section{The Link and Tor} \label{TLAT}

The purpose of this section is to show how it is possible to relate module extensions and links between prime ideals by using homological methods. While the conclusions of the results obtained in this section are already known (see \cite[Chapter 12]{GW}), the proofs are new; indeed, the purpose of this section is to show that the existence of a link can be directly deduced by considering the module theoretic structure of these $\mathrm{Ext}$ modules. 

The following is a trivial result that describes the structure of the $\mathrm{Tor}$ functor in low degrees. While we will only use $\mathrm{Tor_{1}}$ in this section, we record the other information for use later.

\begin{Lemma1} \label{tor}
Let $R$ be a ring, and let $P$ and $Q$ be prime ideals. Then, there are the following isomorphisms: 
\begin{itemize}

\item[a.)] $\mathrm{Tor}_{0}^{R}(R/P,R/Q)\cong R/(P+Q).$
\item[b.)] $\mathrm{Tor}_{1}^{R}(R/P,R/Q)\cong {P\cap Q}/PQ.$
\item[c.)] $\mathrm{Tor}_{2}^{R}(R/P,R/Q)\cong \ker(P\otimes_{R} Q \rightarrow PQ).$

\end{itemize}
\begin{proof}
The first statement is a trivial manipulation of the tensor product. For the second one, we have an exact sequence $0\rightarrow Q\rightarrow R\rightarrow R/Q \rightarrow 0$, leading to an exact sequence $0\rightarrow \mathrm{Tor}^{R}_{1}(R/P,R/Q)\rightarrow Q/PQ \rightarrow R/P \rightarrow R/(P+Q) \rightarrow 0$. The statement follows. For the third, we can dimension shift to get that $\mathrm{Tor}_{2}^{R}(R/P,R/Q)\cong \mathrm{Tor}_{1}^{R}(R/P,Q)$. Considering the long exact sequence that arises from the short exact sequence $0\rightarrow P\rightarrow R\rightarrow R/P\rightarrow 0$ now gives us our result.
\end{proof}
\end{Lemma1} 

The key point of Jategaonkar's work is that certain well behaved module extensions allow us to deduce the existence of links between prime ideals. We show that it is possible to express some of these results in the language of homological algebra. It appears that the relationship between extensions of modules with suitable associated primes and the structure of the linking bimodule between these prime ideals comes arises from a spectral sequence.

Suppose $R$, $S$, $U$, and $V$ are rings, and we have bimodules $_{U}L_{R}$, $_{R}M_{S}$ and $_{V}N_{S}$. The standard adjunction between $\mathrm{\otimes}$ and $\mathrm{Hom}$ gives us the following isomorphism of $(V,U)$-bimodules: $$\mathrm{Hom}_{R}(L,\mathrm{Hom}_{S}(M,N))\cong \mathrm{Hom}_{S}(L\otimes_{R} M, N).$$  We can derive this to get the following spectral sequence. 

\begin{Proposition1} \label{ssp}
Let $R$, $S$, $U$, and $V$ be rings, and suppose we have bimodules  $_{U}L_{R}$, $_{R}M_{S}$ and $_{V}N_{S}$. Then there exist two spectral sequences of $(V,U)$-bimodules, $\{\mathrm{^{I}E^{p,q}_{r}}\}_{r\geq 2}$ and  $\{\mathrm{^{II}E^{p,q}_{r}}\}_{r\geq 2}$, converging to the same graded object $H^{n}$, with 
\begin{align*}
\mathrm{^{I}E^{p,q}_{2}\cong Ext}^{p}_{R}(L,\mathrm{Ext}^{q}_{S}(M,N)), \\
\mathrm{^{II}E^{p,q}_{2}\cong Ext}^{p}_{S}(\mathrm{Tor}_{q}^{R}(L,M),N).
\end{align*}

\begin{proof}
There is a standard argument that is used to prove results of this type - see \cite[Chapter 5.6]{WIEB}. However, the subtlety in our situation arises because we wish to construct spectral sequences of bimodules; in general, we do not know whether it is possible to extract resolutions of $L$ and $N$ which are both composed of bimodules and are acyclic with respect to the functors that we are deriving.

To get around this, we note the following fact: the taking of a projective resolution is a functor from the category of right $R$-modules to $\mathrm{K}(R)$, the homotopy category of (co)chain complexes of right $R$-modules (\cite[Theorem 2.26, Theorem 2.37]{WIEB}). Thus, if $L$ is a $(U,R)$-bimodule, and $\mathrm{P}$ is a projective resolution (of $R$-modules) of $L$, then there is a ring homomorphism $U\rightarrow \mathrm{End}_{\mathrm{K}(R)}(P)$. In a similar fashion, let $I$ be an injective resolution of $N$ as an $R$-module. If $N$ is a $(V,S)$-bimodule, then there is a ring homomorphism $V\rightarrow \mathrm{End}_{\mathrm{K}(R)}(I)$. These actions are preserved by the tensor product, taking $\mathrm{Hom}$, and taking homology groups. If both $L$ and $N$ are bimodules, the reader can check that these actions are compatible with each other. Now, we form the standard double complex, and extract the two spectral sequences that emerge, using the adjunction that we described before this result. Let us describe the construction of the first spectral sequence given above in a bit more detail; the construction of the second is similar. 

Form the double complex $C_{p,q}=\mathrm{Hom}_{R}(P_{p},\mathrm{Hom}_{S}(M,I^{q})).$ These abelian groups are not necessarily $(V,U)$-bimodules (the rings act only up to homotopy); however, every element $v\in V$ (resp.~$u\in U$) does act as an endomorphism of this abelian group, and given $f\in \mathrm{Hom}_{R}(P_{p},\mathrm{Hom}_{S}(M,I^{q}))$, $(vf)u=v(fu)$. The maps in the double complex preserve the action of both $u$ and $v$, for all $u\in U$ and $v\in V$. If we look at the spectral sequence arising from the filtration of the above double complex by columns, we see that $$\mathrm{^{I}E^{p,q}_{1}\cong Hom}_{R}(P_{p}, \mathrm{Ext}^{q}_{S}(M,N)).$$ These actually are left $V$-modules, as $V$ does act on the abelian group up to homotopy, and there is still an action of every element $u\in U$ on the right. Again, all of the maps in the spectral sequence preserve these actions. It then follows that the groups on the $\mathrm{^{I}E_{2}}$ page are $(V,U)$-bimodules. Matters for the second spectral sequence are similar. 
\end{proof}
\end{Proposition1}

\begin{Rem}
The above result appears for algebras over a field in \cite[Proposition 1.3]{BL}. When working with algebras defined over a field, it is possible to take suitable resolutions of bimodules, and we can thus avoid working in the homotopy category.
\end{Rem}

In particular, for any indecomposable injective $E$ with assassinator $Q$,  there is a spectral sequence $$\mathrm{E_{2}^{p,q}=Ext}^{p}_{R/Q}(\mathrm{Tor}_{q}^{R}(R/P,R/Q),\pt{E})\Rightarrow \mathrm{\mathrm{T}\Gamma}^{p+q}_{P}(\pt{E}).$$ We are especially interested in $\mathrm{\mathrm{T}\Gamma}_{P}^{1}(\pt{E})$. As $\pt{E}$ is an injective $R/Q$-module, $\mathrm{E_{2}^{p,q}}=0$  for all $p\geq 1$. It follows that $\mathrm{E_{2}^{1,0}=0}$, and $\mathrm{E_{2}^{0,1}=E_{\infty}^{0,1}}$. Therefore, there is an isomorphism:

\begin{equation} \label{sss}
\mathrm{Hom}_{R/Q}(\mathrm{Tor}^{R}_{1}(R/P,R/Q),\pt{E})\cong \mathrm{\mathrm{T}\Gamma}_{P}^{1}(\pt{E}).
\end{equation}

\vspace{2mm}

With this isomorphism, we can now recover the connection between module extensions and links. 

\begin{Theorems1} \label{duality1}
Let $R$ be a right noetherian ring, and let $P$ and $Q$ be prime ideals. Let $E$ be an indecomposable injective with assassinator $Q$. If $\mathrm{\mathrm{T}\Gamma}_{P}^{1}(\pt{E})$ has a fully faithful $R/P$-submodule, then there exists an ideal $A$ such that $PQ\leq A \lneq P\cap Q$ with $P\cap Q /A$ fully faithful as a left $R/P$-module and as a right $R/Q$-module. If $R$ is noetherian, then $P\rightsquigarrow Q$. 
\begin{proof}
Using \cref{sss} and \cref{tor}, there exists a fully faithful submodule of $$\mathrm{Hom}_{R/Q}(P\cap Q/PQ,\pt{E}),$$ say, $M$. Let $0\neq f\in M$, and take $A$ to be a maximal ideal containing $PQ$ and contained in $P\cap Q$ such that $A$ is contained in $\ker(fr)$ for some $r\in R/P$ such that $fr\neq 0$. It should be noted that in any bimodule $B$, the existence of a one sided unfaithful submodule is equivalent to the existence of a sub-bimodule unfaithful on one side. 

We first show that $P\cap Q/A$ is fully faithful as a right $R/Q$-module. Let $A'$ be a sub-bimodule of $P\cap Q/A$ which is unfaithful on the right. We have that there exists an ideal $I$ such that $(A'/A)I=0$. As $\pt{E}$ is a fully faithful $R/Q$-module, $A'$ is contained in $\ker(fr)$, contradicting the maximality of $A$. Thus $P\cap Q/A$ is fully faithful as a right $R/Q$-module. 

 We now proceed to show that $P\cap Q/A$ is fully faithful as a left $R/P$-module. Observe that $\mathrm{Hom}_{R/Q}(P\cap Q/A, \pt{E})$ embeds in $\mathrm{Hom}_{R/Q}(P\cap Q/PQ, \pt{E})$. Let $fr=g$. As $A$ is contained in $\ker(g)$, we have that $gR$ is a fully faithful submodule of $\mathrm{Hom}_{R/Q}(P\cap Q/A, \pt{E})$. Suppose $B$ is a sub-bimodule of $P\cap Q/A$ which is unfaithful as an $R/P$-module.  Note that $P\cap Q/A$ is faithful as a left $R/P$-module; if not, we would contradict the faithfulness of $gR$. If $B$ is non-zero, then by the maximality of $A$, it follows that $B$ is not contained in $\ker(gr)$ for any $r$ such that $gr\neq 0$. 

There is an exact sequence $$0\rightarrow \mathrm{Hom}_{R/Q}(P\cap Q/B,\pt{E})\rightarrow \mathrm{Hom}_{R/Q}(P\cap Q/A, \pt{E})\rightarrow \mathrm{Hom}_{R/Q}(B/A,\pt{E})\rightarrow 0.$$
Since $B/A$ is unfaithful as a left $R/P$-module, there exists an $I\gneq P$ such that $IB$ is contained in $A$. This implies that the last term in the above exact sequence is unfaithful as an $R/P$-module, because $I$ annihilates it. The fact that $B$ is not contained in $\ker(gr)$ for any $r$ such that $gr\neq 0$ implies that $gR$ embeds into $\mathrm{Hom}_{R/Q}(B/A,\pt{E})$. However, this contradicts the fact that $gR$ is a faithful $R/P$-module. Thus $B=0$, and $P\cap Q/ A$ is fully faithful as an $R/P$-module.

If $R$ is noetherian, then \cite[Lemma 8.1]{GW} tells us that we have a link $P\rightsquigarrow Q$.
\end{proof}
\end{Theorems1} 

For tame modules, there is a slightly stronger variant of \cref{duality1}.

\begin{Theorems1} \label{duality3}
Let $R$ be a right noetherian ring, and let $P$ and $Q$ be prime ideals. Let $E_{Q}$ be the tame  indecomposable injective with assassinator $Q$. If $\mathrm{\mathrm{T}\Gamma}_{P}^{1}(\pt{E_{Q}})$ has a torsionfree $R/P$-submodule, then $P\rightsquigarrow Q$.
\begin{proof}
Again, we can proceed as in \cref{duality1}. In this case, the fact that $\pt{E_{Q}}$ is torsionfree allows us to show that our bimodule is torsionfree on the right. We then argue as in \cref{duality1} to get fully faithfulness on the left.
\end{proof}
\end{Theorems1}

The above results tell us that in general, nice structural properties of our local cohomology modules allow us to see the existence of links. For tame modules, there is a converse to \cref{duality3}: links are reflected in the local cohomology.

\begin{Theorems1} \label{rightlinks}
Let $R$ be a right noetherian ring, and let $P$ and $Q$ be prime ideals. Let $E_{Q}$ be the indecomposable tame injective with assassinator $Q$. Then, $\mathrm{\mathrm{T}\Gamma}_{P}^{1}(\pt{E_{Q}})$ has a torsionfree $R/P$-submodule if there exists an ideal $A$ such that $PQ\leq A \lneq P\cap Q$ with $P\cap Q /A$ faithful as a left $R/P$-module and torsionfree as a right $R/Q$-module.
\begin{proof}
Suppose there exists an ideal $A$ such that $PQ\leq A \lneq P\cap Q$ with $P\cap Q /A$ faithful as a left $R/P$-module and torsionfree as a right $R/Q$-module. As $P\cap Q/A$ is torsionfree as an $R/Q$-module, we know that there is an injective map $f: P\cap Q/A\rightarrow \pt{E_{Q}}^{n}$ for some $n$. Let this map decompose into $f_{1},...,f_{n}$, with each $f_{i}\in \mathrm{Hom}_{R/Q}(P\cap Q/PQ,\pt{E_{Q}})$. If $\mathrm{Hom}_{R/Q}(P\cap Q/PQ,\pt{E_{Q}})$ is torsion as an $R/P$-module, we would have that there exists $c_{i}\in C_{R}(P)$ such that $f_{i}c_{i}=0$. Taking a common denominator, we see that there exists a $c\in C_{R}(P)$ such that $fc=0$. This implies that $f(cx)=0$ for all $x\in P\cap Q/A$. As $f$ is injective, we have that $cx=0$ for all $x$, i.e. that $P\cap Q/A$ is unfaithful, which is a contradiction. Thus $\mathrm{Hom}_{R/Q}(P\cap Q/PQ,\pt{E_{Q}})$ is not torsion, and as $R/P$ is right noetherian, it follows that there exists a torsionfree submodule. The isomorphism \cref{sss} now completes the argument.
\end{proof}
\end{Theorems1}

We can now fully recover the conclusions of Jategaonkar's main lemma, \cref{jml}. If we are in the second case of the theorem, then $\mathrm{\mathrm{T}\Gamma}_{P}^{1}(\pt{E})$ has a fully faithful $R/P$-submodule, thereby implying that $P\rightsquigarrow Q$ by \cref{duality1}. It appears as if the local cohomology acts as a bridge between module extensions and links. This perspective is useful because it tells us exactly where in the argument we need to assume that our ring is two sided noetherian. The relationship between module extensions and local cohomology works over one sided noetherian rings, while the relationship between local cohomology and the linking bimodule requires a two sided hypothesis.

\section{$\mathrm{\mathrm{T}\Gamma}_{P}^{1}(\pt{E})$ and Incomparability Conditions}

The aim of this section is to convince the reader that the structure of $\mathrm{\mathrm{T}\Gamma}_{P}^{1}(\pt{E_{Q}})$ as an $R/P$-module has a strong correlation with the structure of the clique of $Q$. We begin with a general lemma on bimodules.

\begin{Lemma1} \label{bi}
Let $S$ and $T$ be rings, and suppose $T$ is prime right Goldie. Let $B$ be an $(S,T)$-bimodule which is of finite length as a left $S$-module. Then, the following are equivalent:
\begin{itemize}
\item[i.)] $B$ is fully faithful as a $T$-module.
\item[ii.)] $B$ is torsionfree as a $T$-module.
\item[iii.)] $B$ is divisible as a $T$-module.
\item[iv.)] $B$ is injective as a $T$-module.
\end{itemize}
\begin{proof}
If $B$ is an $(S,T)$-bimodule with finite length as an $S$-module, and $T$ is prime right Goldie, then $B$ is torsionfree if and only if it is divisible, as right multiplication by a regular element of $T$ will be surjective if and only if it is injective as an $S$-module map from $B$ to itself. Thus ii.) and iii.) are equivalent. Over a semiprime right Goldie ring, torsionfree divisible modules are injective (\cite[Proposition 7.11]{GW}), and injective modules are divisible. Thus iv.) is also equivalent to the above statements. Finally, torsionfree modules over prime right Goldie rings are fully faithful (\cite[Lemma 7.22]{GW}), which gives us the implication ii.) $\Rightarrow$ i.), while i.) $\Rightarrow$ ii.) follows from \cite[Lemma 8.3]{GW}.
\end{proof}
\end{Lemma1}

\begin{Coroll1} \label{goodbi}
Let $R$ be a right noetherian ring. Let $P$ and $Q$ be prime ideals, and consider $\mathrm{\mathrm{T}\Gamma}_{P}^{i}(\mathcal{Q}(R/Q))$. Then, $\mathrm{\mathrm{T}\Gamma}_{P}^{i}(\mathcal{Q}(R/Q))$ is a $(\mathcal{Q}(R/Q),R/P)$-bimodule of finite length as a $\mathcal{Q}(R/Q)$-module. Furthermore, the following are equivalent: 
\begin{itemize}
\item [i.)]  $\mathrm{\mathrm{T}\Gamma}_{P}^{i}(\mathcal{Q}(R/Q))$ is fully faithful as an $R/P$-module.
\item [ii.)] $\mathrm{\mathrm{T}\Gamma}_{P}^{i}(\mathcal{Q}(R/Q))$ is torsionfree as an $R/P$-module.
\item [iii.)] $\mathrm{\mathrm{T}\Gamma}_{P}^{i}(\mathcal{Q}(R/Q))$ is divisible as an $R/P$-module.
\item [iv.)] $\mathrm{\mathrm{T}\Gamma}_{P}^{i}(\mathcal{Q}(R/Q))$ is injective as an $R/P$-module.
\end{itemize}

As all of the above properties are preserved under summands and finite direct sums, we also have a similar statement for $\mathrm{\mathrm{T}\Gamma}_{P}^{i}(\pt{E_{Q}})$.

\begin{proof}
It is immediate that $\mathrm{\mathrm{T}\Gamma}_{P}^{i}(\mathcal{Q}(R/Q))$ is a $(\mathcal{Q}(R/Q),R/P)$-bimodule. As $\mathcal{Q}(R/Q)$ is simple artinian, and $R$ is right noetherian, it follows that $\mathrm{\mathrm{T}\Gamma}_{P}^{i}(\mathcal{Q}(R/Q))$ is finitely generated, and thus of finite length as a $\mathcal{Q}(R/Q)$-module. We can now apply \cref{bi}.
\end{proof}
\end{Coroll1}

We now collect a few results that describe how the structure of the local cohomology modules impacts the clique structure of the ring.

\begin{Proposition1} \label{ext2}
Let $R$ be a right noetherian ring, and let $P$ and $Q$ be prime ideals, with $P$ not strictly contained in $Q$. Then,
\begin{itemize}
\item $\mathrm{\mathrm{T}\Gamma}_{P}^{1}(\pt{E_{Q}})= 0$ if and only if there does not exist a prime ideal $P'$ containing $P$ such that $P'\rightsquigarrow Q$.
\item $\mathrm{\mathrm{T}\Gamma}_{P}^{1}(\pt{E_{Q}})\neq 0$ and has torsion as an $R/P$-module if and only if there exists a prime ideal $P'\gneq P$ such that $P'\rightsquigarrow Q$.
\item $\mathrm{\mathrm{T}\Gamma}_{P}^{1}(\pt{E_{Q}})\neq 0$ has a torsionfree submodule as an $R/P$-module if and only if $P\rightsquigarrow Q$.
\end{itemize}
\begin{proof}
As $P$ is not strictly contained in $Q$, \cref{mainlemma} implies that $\mathrm{ann}_{E_{Q}/\pt{E_{Q}}}(P)\cong \mathrm{\mathrm{T}\Gamma}_{P}^{1}(\pt{E_{Q}})$. Thus, for any prime $P'\geq P$, we see that $\mathrm{ann}_{\mathrm{\mathrm{T}\Gamma}_{P}^{1}(\pt{E_{Q}})}(P')\cong \mathrm{\mathrm{T}\Gamma}_{P'}^{1}(\pt{E_{Q}})$. The result now follows from \cref{duality3}, \cref{rightlinks}, and \cref{goodbi}.
\end{proof}
\end{Proposition1}

The lack of an analogue of \cref{rightlinks} for wild modules means that the corresponding statement must necessarily be weaker. 

\begin{Proposition1} \label{ext1}
Let $R$ be a noetherian ring, and let $P$ and $Q$ be prime ideals, with $P$ not strictly contained in $Q$. Let $E$ be an indecomposable injective with assassinator $Q$. 
\begin{itemize}
\item If $\mathrm{\mathrm{T}\Gamma}_{P}^{1}(\pt{E})\neq 0$ and has an unfaithful submodule as an $R/P$-module, then there exists a $P'\gneq P$ such that $P'\rightsquigarrow Q$.
\item If $\mathrm{\mathrm{T}\Gamma}_{P}^{1}(\pt{E})$ has a fully faithful submodule as an $R/P$-module, then $P\rightsquigarrow Q$.
\end{itemize}
\begin{proof}
The argument is similar to \cref{ext2}. We will need to use \cref{duality1}. Indeed, the second statement is \cref{duality1}. It is included here for completeness. 
\end{proof}
\end{Proposition1}

Under certain conditions, there is more regularity in the structure of the local cohomology of wild modules. We will add to the above result later in the paper (see \cref{impact2} and \cref{impact3}).

A pair of prime ideals $P$ and $Q$ are said to be comparable if one is contained in the other. It is a major and long-standing open question in ring theory as to whether a link can exist between a comparable pair of prime ideals, or whether comparable primes can live in the same clique. We say that a ring satisfies $(INC1)$ if there does not exist a link $P\rightsquigarrow Q$, with $P$ properly containing $Q$. We say that a ring satisfies $(INC2)$ if given prime ideals $P\lneq P'$ and a prime ideal $Q$, $P\rightsquigarrow Q$ implies that $P'$ does not link to $Q$. There are no known examples of a noetherian ring which does not satisfy $(INC1)$ or $(INC2)$. 

By a dimension function, we mean one in the sense of \cite[Definition 6.8.4]{MR}. A ring $R$ is said to have a symmetric dimension function if there are dimension functions on the category of finitely generated right and left modules, both denoted by $\partial$, such that for any bimodule $B$ which is finitely generated on both sides, $\partial(_{R}B)=\partial(B_{R})$. If our ring has a symmetric dimension function, then it satisfies both of the above conditions, as linked prime ideals must have the same value under this function. 

A ring which satisfies the second layer condition must satisfy $(INC1)$ and $(INC2)$:  by \cite[Corollary 14.5]{GW}, if $P$ and $Q$ are linked prime ideals, then the classical Krull dimension of $R/P$ coincides with that of $R/Q$.  

We now prove a theorem giving some equivalent conditions for the incomparability of linked prime ideals. We will need the following easy observation:

\begin{Lemma1} \label{bimoco}
Let $M$ be an $(R,S)$-bimodule, and let $c\in R$. Write $c^{*}:M\rightarrow M$ for left multiplication by $c$.  Let $E$ be an injective $S$-module. Then, given an $S$-module homomorphism $f:M\rightarrow E$, there exists an $S$-module homomorphism $g: M\rightarrow E$ such that 
\begin{displaymath}
\xymatrix{
M \ar[r]^{c^{*}} \ar[d]_{f} & M \ar@{.>}[dl]^{g} \\
E
}
\end{displaymath}
if and only if $f(\ker(c^{*}))=0$.
\begin{proof}
The proof is elementary homological algebra. By the injectivity of $E$, we have that such a map $g$ exists if and only if there exists a map completing
\begin{displaymath}
\xymatrix{
M \ar[r]^{c^{*}} \ar[d]_{f} & cM  \\
E
}
\end{displaymath}

There is an exact sequence $$0\rightarrow \mathrm{Hom}_{R}(cM,E)\rightarrow \mathrm{Hom}_{R}(M,E)\rightarrow \mathrm{Hom}_{R}(\ker(c^{*}),E).$$ Thus $f$ lies in the image of the first map if and only if $f$ lies in the kernel of the last map, i.e. if $f(\ker(c^{*}))=0$.
\end{proof}
\end{Lemma1}

The following corollary should be contrasted with \cite[(3.5)``Injective Producing Lemma"]{LAM}. 

\begin{Coroll1} \label{biinj}
Let $M$ be an $(R,S)$-bimodule. Let $E$ be an injective $S$-module. If $M$ is torsionfree as a left $R$-module, then $\mathrm{Hom}_{S}(M,E)$ is divisible as a right $R$-module.
\begin{proof}
Immediate from \cref{bimoco}, as now $\ker(c^{*})=0$ for all $c\in C_{R}(0)$.
\end{proof}
\end{Coroll1}  

We now apply the above results to the situation we are interested in: 

\begin{Theorems1} \label{INC2}
Let $R$ be a right noetherian ring, and let $P$ and $Q$ be prime ideals, with $P$ not strictly contained in $Q$. Then, the following are equivalent:
\begin{itemize}
\item[i.)] There does not exist a prime ideal $P'$ such that $P'\gneq P$ and $P'\rightsquigarrow Q$.
\item[ii.)]  The set of left torsion elements of $P\cap Q/PQ$ is contained in the right torsion bimodule of $P\cap Q/PQ$.
\item[iii.)] $\mathrm{\mathrm{T}\Gamma}_{P}^{1}(\pt{E_{Q}})$ is injective as a right $R/P$-module.
\item[iv.)] $\mathrm{Tor}_{1}^{R/P}(R/I,P\cap Q/PQ)$ is torsion as a right $R/Q$-module for any left ideal $I\leq R/P$.
\item[v.)] $\mathrm{\mathrm{T}\Gamma}_{P'}^{1}(\pt{E_{Q}})=0$ for all $P'\gneq P$.
\item[vi.)] $\mathrm{Ext}^{1}_{R}(\text{---},\mathcal{Q}(R/Q))$ is an exact functor on the category of $R/P$-modules.
\end{itemize}
\begin{proof}

$i.) \iff v.)$

This is the content of \cref{ext2}.

$v.) \iff iii.)$

\cref{goodbi} implies that $\mathrm{\mathrm{T}\Gamma}_{P}^{1}(\pt{E_{Q}})$ is injective if and only if it is fully faithful. The rest of the argument is, as above, a consequence of \cref{ext2}.

$iii.)  \iff iv.)$

Recall the spectral sequence of \cref{ssp}. We use it to construct the following isomorphism: \begin{align*}
  & \mathrm{Ext}^{n}_{R/P}(R/I,\mathrm{Hom}_{R/Q}(P\cap Q/PQ,\pt{E_{Q}})) \\
\cong & \mathrm{Hom}_{R/Q}(\mathrm{Tor}_{n}^{R/P}(R/I,P\cap Q/PQ),\pt{E_{Q}}).
\end{align*}

It is immediate from the above isomorphism that $\mathrm{Hom}_{R/Q}(P\cap Q/PQ,\pt{E_{Q}})$ is injective as an $R/P$-module  if and only if $ \mathrm{Hom}_{R/Q}(\mathrm{Tor}_{1}^{R/P}(R/I,P\cap Q/PQ),\pt{E_{Q}})=0$ for all left ideals $I$ of $R/P$. This happens if and only if $\mathrm{Tor}_{1}^{R/P}(R/I,P\cap Q/PQ)$ is torsion as a right $R/Q$-module.

$ii.) \iff iii.)$
Given $c\in C_{R}(P)$, let $c^{*}$ denote the left multiplication by $c$ on the biomodule $P\cap Q/PQ$.

Recall the isomorphism from \cref{sss}: 

 $$\mathrm{Hom}_{R/Q}(P\cap Q/PQ,\pt{E_{Q}})\cong \mathrm{\mathrm{T}\Gamma}_{P}^{1}(\pt{E_{Q}}).$$

By \cref{bimoco}, we have that $\mathrm{Hom}_{R/Q}(P\cap Q/PQ,\pt{E_{Q}})$ is divisible as an $R/P$-module if and only if, for all $f \in \mathrm{Hom}_{R/Q}(P\cap Q/PQ,\pt{E_{Q}})$, $f(\ker(c^{*}))=0$ for all $c\in C_{R}(P)$. Note that $\bigcup_{c\in C_{R}(P)}\ker(c^{*})$ is just the set of left torsion elements of $P\cap Q/PQ$, which we denote by $_{l}T$. $\mathrm{Hom}_{R/Q}(P\cap Q/PQ,\pt{E_{Q}})$ is divisible if and only if $f(_{l}T)=0$ for all $f\in \mathrm{Hom}_{R/Q}(P\cap Q/PQ,\pt{E_{Q}})$. Noting the fact that $\ker(c^{*})$ is a right submodule of $P\cap Q/PQ$ and the injectivity of $\pt{E_{Q}}$ as an $R/Q$-module, we see that this can happen if and only if $_{l}T$ is contained in the right torsion sub-bimodule of $P\cap Q/PQ$, as every non-zero non-torsion $R/Q$-module admits a non-zero homomorphism to $\pt{E_{Q}}$.

$iii.)  \iff vi.)$

From the spectral sequence described in \cref{ssp}, it follows that for any $R/P$-module $M$, $$\mathrm{Hom}_{R/P}(M,\mathrm{Ext}^{1}_{R}(R/P,\mathcal{Q}(R/Q)))\cong \mathrm{Ext}^{1}_{R}(M,\mathcal{Q}(R/Q)).$$ Thus $\mathrm{Ext}^{1}_{R}(R/P,\mathcal{Q}(R/Q))$ is injective as an $R/P$-module if and only if $\mathrm{Ext}^{1}_{R}(M,\mathcal{Q}(R/Q))$ is exact as a functor on the category of $R/P$-modules.
\end{proof}
\end{Theorems1}

In particular, we get the following corollary, which gives a set of equivalent conditions for the incomparability of linked primes:

\begin{Coroll1} \label{INC}
Let $R$ be a right noetherian ring, and let $Q$ be a prime ideal. Then, the following are equivalent:
\begin{itemize}
\item[i.)] There does not exist a prime ideal $P$ such that $P\gneq Q$ and $P\rightsquigarrow Q$.
\item[ii.)]  The set of left torsion elements of $Q/Q^{2}$ is contained in the right torsion bimodule of $Q/Q^{2}$.
\item[iii.)] $\mathrm{\mathrm{T}\Gamma}_{Q}^{1}(\pt{E_{Q}})$ is injective as a right $R/Q$-module.
\item[iv.)] $\mathrm{Tor}_{1}^{R/Q}(R/I,Q/Q^{2})$ is torsion as a right $R/Q$-module for any left ideal $I\leq R/Q$.
\item[v.)] $\mathrm{\mathrm{T}\Gamma}_{P}^{1}(\pt{E_{Q}})=0$ for all $P\gneq Q$.
\item[vi.)] $\mathrm{Ext}^{1}_{R}(\text{---},\mathcal{Q}(R/Q))$ is an exact functor on the category of $R/Q$-modules.
\end{itemize}
\begin{proof}
This follows by taking $P$ equal to $Q$ in \cref{INC2}.
\end{proof}
\end{Coroll1}

Let $t(M)$ denote the torsion submodule of $M$.

\begin{Coroll1}
Let $R$ be a noetherian ring and let $Q$ be a prime ideal. Then $Q$ does not link to any comparable prime ideal if and only if $t(_{R/Q}(Q/Q^{2}))=t((Q/Q^{2})_{R/Q})$.
\begin{proof}
This follows using left-right symmetry and \cref{INC}.ii. 
\end{proof}
\end{Coroll1} 

The following corollary can probably be proved by other means, but is an immediate consequence of \cref{INC}.

\begin{Coroll1}
 Let $R$ be a right noetherian ring of right global dimension $1$. Let $Q$ be a prime ideal in $R$. Then, there does not exist a prime ideal $P$ such that $P\gneq Q$ and $P\rightsquigarrow Q$. In particular, if $R$ is a noetherian ring of global dimension $1$, linked prime ideals are incomparable.
\begin{proof}
In this case, it is easy to verify condition iii.) in \cref{INC}. Recall, from \cref{mainlemma}, that $\mathrm{\mathrm{T}\Gamma}_{Q}^{1}(\pt{E_{Q}})\cong \mathrm{ann}_{E/\pt{E}}(Q).$ If $R$ has right global dimension $1$, it follows that $E/\pt{E}$ is an injective $R$-module. Thus $\mathrm{ann}_{E/\pt{E}}(Q)$ is an injective $R/Q$-module.
\end{proof}
\end{Coroll1}

We now prove a result that describes how the structure of the tame indecomposable module with assassinator $P$ has some impact on the structure of wild modules with assassinator $P$. 

\begin{Proposition1} \label{impact1}
Let $R$ be a right noetherian ring, and let $P$ and $Q$ be prime ideals, with $P$ not properly contained in $Q$. Let $E$ be an indecomposable injective module with assassinator $Q$. Then, $\mathrm{\mathrm{T}\Gamma}_{P}^{1}(\pt{E_{Q}})$ embeds as a submodule of $\prod_{i\in I} \mathrm{\mathrm{T}\Gamma}_{P}^{1}(\pt{E})$ for some indexing set $I$.
\begin{proof}
Let $x\in E_{Q}$. The submodule $xR$ contains a submodule $U$ which is isomorphic to a uniform right ideal of $R/Q$. To see this explicitly, note that $xR$ must intersect $\pt{E_{Q}}$ non-trivially. Thus, it contains a submodule which is torsionfree as an $R/Q$-module, and thus contains a uniform right ideal $U$ of $R/Q$. By \cite[Lemma 14.9]{GW}, there is a non-zero homomorphism from $U$ to $\pt{E}$, which extends to a homomorphism $f_{x}:E_{Q}\rightarrow E$ such that $f_{x}(x)\neq 0$. The maps $\{f_{x}\}_{x\in E_{Q}}$ induce an embedding of $E_{Q}$ into $\prod_{i\in I} E$ for some indexing set $I$ (indeed, $I$ can be taken to be $E_{Q}$). It is immediate that $\mathrm{ann}_{\prod_{i\in I} E}(Q)= \prod_{i\in I} \pt{E}$, and we therefore have that $E_{Q}/\pt{E_{Q}}$ embeds into $\prod_{i\in I} E/\pt{E}$. This implies that $\mathrm{ann}_{E_{Q}/\pt{E_{Q}}}(P)$ embeds into $\prod_{i\in I} \mathrm{ann}_{E/\pt{E}}(P)$. As $P$ is not contained in $Q$, however, this is precisely saying that $\mathrm{\mathrm{T}\Gamma}_{P}^{1}(\pt{E_{Q}})$ embeds as a submodule of $\prod_{i\in I} \mathrm{\mathrm{T}\Gamma}_{P}^{1}(\pt{E})$. 
\end{proof}
\end{Proposition1}

\begin{Coroll1} \label{impact2}
Let $R$ be a right noetherian ring. If $\mathrm{\mathrm{T}\Gamma}_{P}^{1}(\pt{E_{Q}})\neq 0$, then $\mathrm{\mathrm{T}\Gamma}_{P}^{1}(\pt{E})\neq 0$ for all indecomposable injective modules $E$ with assassinator $Q$.
\begin{proof}
Immediate from \cref{impact1}
\end{proof}
\end{Coroll1}

Under sufficiently nice conditions, we can say even more about the structure of local cohomology modules of wild modules.

\begin{Proposition1} \label{impact3}
Let $R$ be a noetherian ring which satisfies (INC1) and (INC2). Then, if $P\rightsquigarrow Q$, $\mathrm{\mathrm{T}\Gamma}_{P}^{1}(\pt{E})$ is a fully faithful divisible $R/P$-module for any indecomposable injective $E$ with assassinator $Q$.
\begin{proof}
If $P\rightsquigarrow Q$, then by (INC1), $P$ is not properly contained in $Q$ and $\mathrm{\mathrm{T}\Gamma}_{P}^{1}(\pt{E})\neq 0$, by \cref{impact2}. If $\mathrm{\mathrm{T}\Gamma}_{P}^{1}(\pt{E})$ is not fully faithful as an $R/P$-module, then it must have an associated prime $P'$ which properly contains $P$. As in \cref{INC}, this implies that $\mathrm{\mathrm{T}\Gamma}_{P'}^{1}(\pt{E})$ is a fully faithful $R/P'$-module. However, this implies that there exists a link $P'\rightsquigarrow Q$, by \cref{ext1}. This contradicts (INC2).  Thus $\mathrm{\mathrm{T}\Gamma}_{P}^{1}(\pt{E})$ is a fully faithful $R/P$-module. 

We now show divisibility. By what we have done above, in particular, $\mathrm{\mathrm{T}\Gamma}_{P}^{1}(\pt{E_{Q}})$ is fully faithful, and \cref{goodbi} tells us that it is divisible. If we let $A$ be the right torsion bimodule of $P\cap Q$, then an argument similar to that used in \cref{INC2} tells us that  $P\cap Q/ A$ is torsionfree on both sides. As $R$ is noetherian, $A$ is unfaithful on both sides, and we have, for all indecomposable injectives $E$ with assassinator $Q$, an isomorphism $\mathrm{Hom}_{R}(P\cap Q/ PQ,E)\cong \mathrm{Hom}_{R}(P\cap Q/ A, E)$. It follows from \cref{biinj} that $\mathrm{Hom}_{R}(P\cap Q/ PQ,E)$ is a divisible $R/P$-module. The spectral sequence in \cref{ssp} tells us that this is isomorphic to $\mathrm{\mathrm{T}\Gamma}_{P}^{1}(\pt{E})$.
\end{proof}
\end{Proposition1}

\section{Links between Small Prime Ideals}

In this section, we aim to use some of results in this paper to study the link structure of primes ideals of small homological height. This question was studied for height 1 prime ideals in various rings satisfying a polynomial identity by Braun in \cite{BRA1}, Braun and Warfield in \cite{BRAWAR}, and Lenagen in \cite{LEN1}. Our aim will be to prove results regarding the local link structure of localizable cliques in homologically well behaved noetherian rings. 

We will need to assume that the reader is familiar with the theory of Auslander-Gorenstein rings, and the notion of an AS-Gorenstein ring. As a reference for details about these objects, we recommend \cite{z-f}.

Over a right noetherian ring, there is an integer valued function $j_{R}$ on the category of finitely generated modules, known as the \textit{grade}. $$j_{R}(M) = \mathrm{inf}\{i|\mathrm{Ext}^{i}_{R}(M,R)\neq 0\}.$$

Over an Auslander-Gorenstein ring $R$, $-j_{R}$ is a dimension function, in the sense of \cite[Definition 6.8.4]{MR}. An Auslander-Gorenstein ring $R$ is said to be \textit{grade-symmetric} if given any $(R,R)$-bimodule $B$ which is finitely generated on both sides, $j_{R}(_{R}B)=j_{R}(B_{R})$. $R$ is said to be \textit{weakly bifinite} if given a subquotient $B$ of $R$, $\mathrm{Ext}^{i}_{R}(B,R)$ is finitely generated as a module on both sides for all $i$.

We record one more definition, for use later in this paper. An artinian ring $A$ is said to be \textit{weakly symmetric} if given rings $L$ and $T$ and bimodule $_{L}M_{A}$ and $_{T}N_{A}$ (resp.~$_{A}M_{L}$ and $_{A}N_{T}$) which are finite length on both sides, $\mathrm{Hom}_{A}(M,N)$ is a $(T,L)$ (resp.~$(L,T)$) bimodule which is finite length on both sides.

If $M$ is a module over a ring $R$, the \textit{socle} series (\cite[Page 78]{GW}) of $M$ is the ascending chain of submodules defined inductively by $\mathrm{soc}^{n+1}(M)/\mathrm{soc}^{n}(M)=\mathrm{soc}(M/\mathrm{soc}^{n}(M))$.

\begin{Definitions1}
 Let $\mathrm{C}$ be a clique of prime ideals in a noetherian ring $R$. We say that $\mathrm{C}$ is \textit{localizable} if $X$, the intersection of elements regular modulo $Q$ for every prime $Q\in \mathrm{C}$, is a right and left Ore set, and the localization $R_{X}$ satisfies the following properties:
\begin{itemize}
 \item[a.)] Prime ideals of the form $Q_{X}$, as $Q$ varies across $\mathrm{C}$, are the only right and left primitive ideals in $R_{X}$.
 \item[b.)] The quotient ring $R_{X}/Q_{X}$ is simple artinian for all $Q\in C$.
\end{itemize}
$X$ is said to be \textit{right classically localizable} if, in addition to the conditions given above, for all $P\in X$, $E_{P_{X}}$ is the union of its socle series. If $C$ is right (classically) localizable, with corresponding Ore set $X$, we will denote its localization by $R_{X}$.
\end{Definitions1}

Suppose we have a localizable clique $C$ in a ring $R$; let its corresponding Ore set be denoted by $X$. In this case, as the only left and right primitive ideals of $R_{X}$ are those of the form $Q_{X}$ for $Q\in C$, the left and right simple $R_{X}$-modules are parametrized by the prime ideals in the clique. For the purposes of this section, we will use the notation $S^{P}$ to refer to the unique simple right $R_{X}$-module with annihilator $P_{X}$, and ${^{P}S}$ for the unique simple left $R_{X}$-module with the same property.

In certain situations, the localization of a ring at a localizable clique has nice homological properties.

\begin{Proposition1} \label{locduality}
Let $R$ be an Auslander-Gorenstein, grade-symmetric ring, and let $C$ be a localizable clique of prime ideals, with corresponding Ore set $X$.  Then, $R_{X}$ is an AS-Gorenstein ring of injective dimension $j_{R}(R/P)$, where $P$ is any prime ideal in $C$, and there is a bijection $D:C\rightarrow C$ such that given $P\in C$ and any right module $M$ of finite flat dimension, there is an isomorphism: $$\mathrm{Tor}_{i}^{R_{X}}(M_{X}, {^{D(P)}S)}\cong \mathrm{Ext}^{j_{R}(R/P)-i}_{R_{X}}(S^{P},M_{X}).$$
\begin{proof}
This is \cite[Proposition 1.10]{VYAS1}. Another reference is \cite[Proposition 2.2.2]{VYASTHE}.
\end{proof}
\end{Proposition1}

We note that there is a symmetric result for left modules, with the map $D$ replaced by the map $D'$.

The following lemma is also useful; it tells us that links are reflected in localizations.

\begin{Lemma1} \cite[Exercise 13H]{GW} \label{local5}
Let $X$ be a right Ore set in a noetherian ring. Let $P$ and $Q$ be prime ideals disjoint from $X$. Then, $P_{X}\rightsquigarrow Q_{X}$ if and only if $P\rightsquigarrow Q$.
\end{Lemma1}

The following is an immediate consequence of what we have done above. Note that if $Q$ belongs to a localizable clique with Ore set $X$, then $\pt{E_{Q}}$ carries an $R_{X}$-module structure, and is isomorphic to a simple summand of $R_{X}/Q_{X}$, i.e. it is isomorphic to $S^{Q}$.

\begin{Proposition1} \label{smalllinks}
Let $R$ be a grade-symmetric, Auslander regular ring. Let $C$ be a localizable clique, and let $P$ and $Q$ be prime ideals in $C$. Then, $$\mathrm{Ext}^{1}_{R_{X}}(R_{X}/P_{X}, S^{Q})\cong \mathrm{Tor}^{R_{X}}_{j_{R}(R/P)-1}(S^{Q}, {^{D(P)}S})^{l(\mathcal{Q}(R/P))}$$ In particular, $P\rightsquigarrow Q$ is and only if $ \mathrm{Tor}^{R_{X}}_{j_{R}(R/P)-1}(R_{X}/Q_{X}, R_{X}/D(P)_{X})\neq 0.$
\begin{proof}
This is immediate from \cref{locduality}, \cref{ext2} and \cref{local5}.
\end{proof}
\end{Proposition1}

\begin{Rem}
 It is sometimes possible to better understand the bijection $D$ mentioned in \cref{locduality}. Whenever $R$ is an dualizing complex over itself with the trace property for prime ideals with respect to a bijection $\phi: \operatorname{Spec}(R)\rightarrow \operatorname{Spec}(R)$, it follows that $D$ is induced by this map $\phi$. This situation often arises when $R$ is a complete semi-local noetherian ring which admits a Morita duality: in this case the trace map is $\phi= \psi\eta$, where $\psi$ is the automorphism of $\operatorname{Spec}(R)$ coming from the Morita duality, and $\eta$ is a distinguished automorphism of $R$ - see \cite{CWZ} and \cite{JY}. Other examples of this phenomena arise when $R$ is a connected graded $k$-algebra over some field - again, \cite{JY} is a good reference for  more information on this situation. 
\end{Rem}

We can use \cref{smalllinks} and \cref{local5} to determine the link structure around primes of small homological height. Recall the isomorphisms given in \cref{tor}.

\begin{Proposition1} \label{height1}
Let $R$ be a grade-symmetric Auslander regular ring. Let $P$ be a prime ideal with a localizable clique such that $j_{R}(R/P)=1$. Then, $P\rightsquigarrow Q$ if and only if $Q=D(P)$. The clique of $P$ looks like: $$\{ ...\rightsquigarrow D'(P)\rightsquigarrow P\rightsquigarrow D(P)\rightsquigarrow... \}.$$
\begin{proof}
By \cref{smalllinks}, we are looking for primes $Q$ such that $$\mathrm{Tor}^{R_{X}}_{0}(R_{X}/Q_{X}, R_{X}/D(P)_{X})\neq 0.$$ Note that in $R_{X}$, the primes in the clique of $P$ correspond to distinct maximal ideals. Thus, $ \mathrm{Tor}^{R_{X}}_{0}(R_{X}/Q_{X}, R_{X}/D(P)_{X})\neq 0$ if and only if $Q=D(P)$. Thus $P\rightsquigarrow Q$ if and only if $Q= D(P)$. The fact that $D'D$ is the identity tells us that $Q\rightsquigarrow P$ if and only if $P=D(Q)$ if and only if $D'(P)=Q$. The result on the structure of the clique follows.
\end{proof}
\end{Proposition1}

In particular, a prime ideal of homological height $1$ is localizable if and only if $D(P)=P$. 

Thus for nice rings, cliques of small prime ideals appear in a `line'. It is curious to note that the situation for prime ideals of homological height $2$ has quite a different flavour to it.

\begin{Proposition1}\label{height2}
Let $R$ be a grade-symmetric Auslander regular ring. Let $P$ be a prime ideal with a localizable clique such that $j_{R}(R/P)=2$. Then, $P\rightsquigarrow Q$ if and only if $Q\rightsquigarrow D(P)$. 
\begin{proof}
In this case, we are looking for prime ideals $Q$ such that $$\mathrm{Tor}^{R_{X}}_{1}(R_{X}/Q_{X}, R_{X}/D(P)_{X})\neq 0.$$ However, as before, the primes in the clique of $P$ correspond to maximal ideals in $R_{X}$. By \cref{tor}, $\mathrm{Tor}^{R_{X}}_{1}(R_{X}/Q_{X}, R_{X}/D(P)_{X})\neq 0$ if and only if $Q_{X}\rightsquigarrow D(P)_{X}$. However, by \cref{local5}, this happens if and only if $Q\rightsquigarrow D(P)$.
\end{proof}
\end{Proposition1}

\section{Density and Indecomposable Injectives}

An immediate consequence of Jategaonkar's main lemma (\cref{jml}) is that over a noetherian ring, the second layer of a tame module contains only tame submodules; we generalized this fact to one sided noetherian rings in \cref{mainlemmafull}. 

In this section, we attempt to study such questions for wild modules. The results that we will obtain give, in particular, some information about the structure of wild indecomposable injective modules. 

\begin{Ques}
Let $R$ be a noetherian ring. Let $U$ be a wild uniform module over $R$. Can the second layer of $U$ contain a tame submodule?
\end{Ques}

By using \cref{mainlemma}, we see that this question is equivalent to the following one:

\begin{Ques}
Let $R$ be a noetherian ring, and let $P$ be a prime ideal in $R$. Is it possible for $\mathrm{\mathrm{T}\Gamma}_{P}^{1}(\pt{E})$ to have a torsionfree submodule when $E$ is a wild indecomposable injective module?
\end{Ques}

Such questions were investigated by Jategaonkar in \cite{JAT2} for polynormal rings, Brown in \cite{KB2} for group rings of polycyclic-by-finite groups, and Lenagan in \cite{LEN2} for algebras with finite Gelfand-Kirilov dimension which satisfy the strong  second layer condition. In all of these cases, it was shown that second layers of wild modules are indeed wild: the proof for polynormal rings appears in \cite[Theorem 3.7]{JAT2}, while the corresponding argument for polycyclic-by-finite group rings, as described in \cite[Theorem 5.12]{KW}, follows by looking at \cite[Corollary 7.2]{KB2} and \cite[Proposition 6.3.14]{JAT}. The result for algebras with finite GK-dimension appears as \cite[Corollary 1.8]{LEN2}.

We also note that there is an example of a ring such that the second layers of wild modules are not wild, constructed by Goodearl and Schofield in \cite{GS}. Both this negative result and the positive results given above were framed in terms of Jategaonkar's \textit{density condition}. Let $R$ and $S$ be prime noetherian rings, and let $B$ be an $(S,R)$-bimodule which is torsionfree on both sides. $B$ is said to satisfy the \textit{right density condition with respect to an indecomposable injective $E$} if $\mathrm{Hom_{R}(B,E)}$ is a torsion $S$-module. We further say that $B$ satisfies the density condition if it satisfies the density condition with respect to every indecomposable injective module. By abuse of notation, a noetherian ring $R$ is said to satisfy the density condition if given prime ideals $P$ and $Q$, every linking bimodule between $P$ and $Q$ satisfies the right density condition. Jategaonkar proved, in \cite[Theorem 3.5]{JAT2}, that a ring satisfies the right density condition if and only if second layers of wild modules are wild. While we do not give a proof of this result in this section, the reader should be convinced of its plausibility by examining the formula \cref{sss}.

The following concepts were defined in \cite{VYAS1} and \cite{VYASTHE}.

\begin{Definitions1}
Let $R$ be a grade-symmetric Auslander-Gorenstein ring. A prime ideal $P$ is \textit{right correct} if $\mathrm{Ext}^{i}_{R}(R/P,R)$ is torsion as an $R/P$-module whenever $i\gneq j_{R}(R/P)$.
\end{Definitions1}

\begin{Definitions1} \label{defofloccor}
Let $R$ be a grade-symmetric Auslander-Gorenstein ring. A prime ideal $P$ is \textit{right localizably correct} if
\begin{itemize}
\item[1.)] The clique of $P$ is localizable, with corresponding Ore set $X$.
\item[2.)] $P$ is right correct.
\item[3.)] There are isomorphisms of $(R_{X},R_{X})$-bimodules: 
\begin{align*}
\mathrm{Ext}^{j_{R}(R/P)}_{R}(R/P,R)\otimes_{R} R_{X} & \cong R_{X}\otimes_{R} \mathrm{Ext}^{j_{R}(R/P)}_{R}(R/P,R) \\
& \cong   \mathrm{Ext}^{j_{R}(R/P)}_{R_{X}}(R_{X}/P_{X},R_{X}).\\
\end{align*}
\end{itemize}
A clique of prime ideals $C$ is \textit{right localizably correct} if all of the prime ideals in $C$ are right localizably correct.
\end{Definitions1}

There are obvious definitions of \textit{left correct} and \textit{left localizably correct} prime ideals, where we instead choose to work in the category of left modules. A prime ideal is \textit{correct} (resp.~\textit{localizably correct}) if it is both left and right correct (resp.~localizably correct). There is the obvious definition of a localizably correct clique.

\begin{Proposition1} \label{listofloc}
The following are examples of correct and localizably correct prime ideals:
 \begin{itemize}
\item [a.)] If $R$ is an Auslander-Gorenstein, grade-symmetric, weakly bifinite ring, then every prime ideal is correct. Furthermore, if $R$ satisfies the second layer condition, then every prime ideal with a localizable clique is localizably correct.
\item [b.)] Let $k$ be a field. Let $R$ be a $k$ algebra with a noetherian connected filtration such that $gr(A)$ is commutative of finite global dimension. Then, every prime ideal in $R$ is correct.  If a prime ideal $P$ has a localizable clique, then it is localizably correct.
\item [c.)] Let $k$ be a field. Let $R$ be a complete local Auslander-Gorenstein $k$-algebra with maximal ideal $\mathfrak{m}$ such that $R/\mathfrak{m}$ is weakly symmetric. Then, every prime ideal in $R$ is correct.  If a prime ideal $P$ has a localizable clique, then it is localizably correct.

\end{itemize}
\begin{proof}
This follows from \cite[Proposition 4.1]{VYAS1} (or \cite[Proposition 3.2.3]{VYASTHE}), using \cite[Theorem 6.17]{JY2}, \cite[Proposition 6.18]{JY2} and \cite[Proposition 4.5 (1)]{JY} for b.), and \cite[Theorem 0.1]{CWZ} and \cite[Proposition 4.5 (2)]{JY} for c.). Part a.) is a direct consequence of \cite[Proposition 4.1]{VYAS1} (or \cite[Proposition 3.2.3]{VYASTHE}).
\end{proof}
\end{Proposition1}

\begin{Rem}
It turns out that the statement in part c.) of the above proposition continues to hold when $k$ is a discrete valuation ring. We show this in a forthcoming paper, tentatively titled `Dualizing Complexes and Discrete Valuation Rings'. Another reference for this is \cite{VYASTHE}.
\end{Rem}

The above notions were defined in order to prove the following theorem.

\begin{Theorems1} \label{localizeext}
Let $R$ be a grade-symmetric Auslander-Gorenstein ring, and let $P$ be a right localizably correct prime ideal with clique $C$ with corresponding Ore set $X$. Let $M$ be a module of finite flat dimension. Then, for all $i$ there is an isomorphism of right $R/P$-modules: $$\mathrm{Ext}^{i}_{R}(R/P,M)\otimes_{R} R_{X}\cong \mathrm{Ext}^{i}_{R_{X}}(R_{X}/P_{X},M_{X}).$$
\begin{proof}
This is \cite[Theorem 4.3]{VYAS1}. Another reference is \cite[Theorem 3.2.5]{VYASTHE}.
\end{proof}
\end{Theorems1}

\begin{Theorems1} \label{rdc2}
Let $R$ be a grade-symmetric Auslander regular ring. Let $P$ be a localizably correct prime ideal. Then, given a wild indecomposable injective module $E$ with assassinator $Q$, $\mathrm{T\Gamma}_{P}^{1}(\pt{E})$ is a torsion $R/P$-module. 
\begin{proof}
Suppose not. If $\mathrm{T\Gamma}_{P}^{1}(\pt{E})$ has a torsionfree submodule, then \cref{ext2} implies that $P\rightsquigarrow Q$. Thus $Q$ lies in the clique of $P$. It follows that $\pt{E}$ vanishes when we localize at the clique of $P$, and by \cref{localizeext} we have a contradiction.
\end{proof}
\end{Theorems1}

\begin{Coroll1}
Let $R$ be a grade-symmetric Auslander regular ring. Let $M$ be a uniform $R$-module with associated prime $Q$ such that $\operatorname{pt}(M)$ is torsion as an $R/Q$-module. Then, given a localizably correct prime $P$, every submodule of $M/\operatorname{pt}(M)$ with associated prime and annihilated by $P$ is torsion as an $R/P$-module.
\end{Coroll1}

We can reformulate \cref{rdc2} into a result involving the density condition. Note that by \cref{bimoco}, in a ring which satisfies (INC1) and (INC2), the left and right torsion sub-bimodules of $P\cap Q/ PQ$ coincide whenever $P\rightsquigarrow Q$. We call this common bimodule the \textit{torsion bimodule} of the link. The quotient of $P\cap Q$ by the torsion bimodule coincides with what Jategaonkar refers to as the `strongest link' between $P$ and $Q$.

\begin{Theorems1} \label{rdc3}
Let $R$ be a grade-symmetric, Auslander regular ring such that every prime ideal is localizably correct. Let $P$ and $Q$ be prime ideals such that $P\rightsquigarrow Q$, and let the clique of $P$ be localizable. Then, $P\cap Q / A$ satisfies the right density condition, where $A$ is the torsion bimodule of $P\cap Q/ PQ$. Furthermore, every linking bimodule between $P$ and $Q$ satisfies the right density condition, and $R$ satisfies the left and right density conditions.
\begin{proof}
There is a bijection between injective modules over $R/Q$ and injective modules over $R$ with assassinators containing $Q$. Thus, let $E$ be an injective module with assassinator containing $Q$, excluding  $E_{Q}$. We have, by the spectral sequence in \cref{ssp}, an isomorphism $$\mathrm{Hom}_{R/Q}(P\cap Q/ A,\mathrm{ann}_{E}(Q))\cong \mathrm{T\Gamma}_{P}^{1}(\mathrm{ann}_{E}(Q)).$$ We can now argue as in \cref{rdc2} to see that the above module is torsion as an $R/P$-module. 

If $P\cap Q/ I$ is any linking bimodule, it is immediate that $A\subset I$. It therefore follows that $\mathrm{Hom}_{R/Q}(P\cap Q/ I,\mathrm{ann}_{E}(Q))$ is a submodule of $\mathrm{Hom}_{R/Q}(P\cap Q/ A,\mathrm{ann}_{E}(Q))$. Thus $P\cap Q/ I$ also satisfies the right density condition.
\end{proof}
\end{Theorems1}

\subsection{Transfer Across Links}
One can ask what properties of quotient rings are transferred across links. For example, \cite[Exercise 13D]{GW} tells us that if $P\rightsquigarrow Q$ in a noetherian ring $R$, then $R/P$ is artinian if and only if $R/Q$ is artinian. The proof is an application of Lenagan's Theorem. We prove a result of a similar flavour:

\begin{Proposition1}
Let $R$ be a noetherian ring satisfying (INC1) and (INC2) with the following property: if $P\rightsquigarrow Q$, and $E$ is a wild indecomposable injective with assassinator $Q$, then $\mathrm{T\Gamma}_{P}^{1}(\pt{E})$ is not torsionfree as an $R/P$-module. If $P\rightsquigarrow Q$, then $R/P$ is right bounded if and only if $R/Q$ is right bounded. 
\begin{proof}
The fact that $R/Q$ right bounded implies $R/P$ right bounded was shown in \cite[Remark 3.6]{KW}. We prove the converse.

Suppose $R/Q$ is not right bounded. Then, by the proof in \cite[Theorem 9.15]{GW}, we see that there exists a wild indecomposable injective $E$ with assassinator $Q$. By our hypothesis, $\mathrm{T\Gamma}_{P}^{1}(\pt{E})$ contains torsion, and by \cref{impact3}, our hypothesis implies that it is fully faithful as an $R/P$-module. These facts together imply that $R/P$ is not right bounded. 
\end{proof}
\end{Proposition1}

The above proposition is, in particular, applicable to rings satisfying the density condition on linking bimodules, such as those described above. 

\section*{Acknowledgments}
This work was done as part of the author's PhD thesis, which was supervised by Simon Wadsley. The author would like to thank Simon Wadsley for his many helpful comments and suggestions while reading through many preliminary versions of this paper. The author would also like to thank Jonathan Nelson for his careful proof-reading of a preliminary version of this paper, and his many suggestions and helpful comments. The author thanks the Cambridge Commonwealth Trust and Wolfson College, Cambridge for their support over the years.

\end{document}